\documentclass[12pt,a4paper]{article}
\usepackage{amssymb, amsmath, amsfonts, amsthm, amscd}
\usepackage[cp1251]{inputenc}
\usepackage[russian]{babel}

\theoremstyle{definition}

\newcommand{\wh}[1]{\widehat{#1}}

\makeatletter
\renewcommand{\@begintheorem}[2]{\par\it {\bf #1\ #2. }}
\renewcommand{\@endtheorem}{}

\renewcommand{\@biblabel}[1]{#1.\hfill}
\makeatother

\def\T{{\mathbb T}}
\def\Ind{\operatorname{Ind}}

\def\Aut{\operatorname{Aut}}

\begin{document}

\begin{center}
\bf
 Twisted Burnside Theory for the Discrete Heisenberg
Group and for Wreath Products of Some Groups
 \end{center}

\begin{center}
{\bf  F.K. Indukaev}\\
Moscow State University\\
\tt gecko984@gmail.com
  \end{center}

\begin{center}{\bf INTRODUCTION.} \\   \end{center} Let $G$ be a countable discrete group and $\phi$
--- its automorphism. Consider the following equivalence
relation on the group $G$.

 {\bf Definition. } Two elements $x,x'\in G$ are said to be
$\phi$-{\em conjugate} or {\em twisted-conjugate,} if there exists
an element $g \in G$ such, that $ x'=g x \phi(g^{-1})$. The set of
elements being $\phi$-conjugate to an element $x$, is called the
$\phi$-{\em conjugacy class} or {\em twisted conjugacy class} class
of the element $x$ and is denoted by $\{x\}_\phi$.

It is easy to see that if $\phi$ is the identity mapping, then the
$\phi$-conjugacy classes are the usual conjugacy classes of the
group $G$.

We are interested in the number of $\phi$-conjugacy classes, which
the group G is partitioned in. It is called the {\em Reidemeister
number} of the automorphism $\phi$ and is denoted by $R(\phi)$. For
different groups and their automorphisms this number can be both
finite and infinite. For example, it was shown in [1] that the
Reidemeister number of any automorphism of a Gromov hyperbolic group
is infinite, and in this paper for the discrete Heisenberg group we
construct automorphisms having any preassigned even Reidemeister
number. Other results concerning Reidemeister numbers may be found
in [2–4].

Consider the set $\wh G_f$ of unitary equivalence classes of
finite-dimensional irreducible unitary representations of the group
$G$. Each automorphism  $\phi$ generates the mapping $\wh\phi_f
\colon \wh G_f \to \wh G_f$, given by the formula $\rho \mapsto \rho
\circ \phi$.

Let $S_f(\phi)$ be the number of fixed points of the mapping
$\wh\phi_f$. In [5], the class RP of groups was defined, for which
this number is equal to $R(\phi)$ providing the latter number is
finite.

This statement is a generalization of the classic Burnside theorem
[6] stating that for a finite group the number of irreducible
representations classes coincides with the number of conjugacy
classes of the group. Indeed, let a group $G$ be finite, and $\phi$
be the identity mapping. Then all representations of the group are
unitary, finite-dimensional, and they remain fixed under the mapping
$\wh\phi_f$, and $R(\phi)$ is the number of usual conjugacy classes.
At last, the class RP contains, in particular, all finite groups
(see [5]).

In former papers [7] and [8] a twisted Burnside type theorem was
proved for groups of type I and for two-step torsion-free nilpotent
groups, respectively. In [9], it was shown that consideration of
only finite-dimensional representations is essential.

It is known that all finite, all Abelian, and all almost polycyclic
groups belong to the class RP (see [5]). But for now there is no
simple criterion for a group to belong to this class. Therefore, it
is reasonable to look for new examples of groups from the class RP
(or, simply, RP-groups). In Section 3 of this paper the RP-property
is proved for wreath products of finitely generated Abelian groups
with the group of integers. Such wreath products become the first
known example of finitely generated RP-groups being not almost
polycyclic.

In the first part of the paper (Sections 1 and 2), the discrete
Heisenberg group is considered. It is shown that for this group all
even numbers can be realized as Reidemeister numbers and also fixed
representations are found explicitly for a specific automorphism
with the Reidemeister number $R(\phi)=2$. \\

\begin{center}{\bf 1. AN INFINITE SEQUENCE OF DISCRETE HEISENBERG GROUP
AUTOMORPHISMS}\\   \end{center}
By H we denote the {\em discrete Heisenberg group}, which is defined
as the following semidirect product:

$$
\begin{aligned}
H=\mathbb{Z}^2 \rtimes \mathbb{Z}, \qquad \widetilde{\alpha} \colon
\mathbb{Z} \to \Aut(\mathbb Z^2); s \mapsto {\alpha}^s, \quad
\alpha=(*,*) \bigl(\begin{smallmatrix} 1 & 1\\ 0&1 \end{smallmatrix}
\bigr).
\end{aligned}
$$
Thus, this group consists of triplets $((m,k),s)$ of integers with
the following multiplication law:
$$
((m,k),s)*((m',k'),s')=((m,k)+\alpha^s(m',k'))=((m+m',k+k'+sm'),s+s');
$$
In particular, $ ((m,k),0)*((0,0),s)=((m,k),s) $. The inverse
element may be found by the formula $ ((m,k),s)^{-1}=((-m,sm-k),-s).
$

The group $H$ can be also considered as the group of integer
$3\times 3$-matrices of the form
$$
\begin{pmatrix}
1 & s & k\\
0 & 1 & m\\
0 & 0 & 1
\end{pmatrix}
$$
with respect to matrix multiplication. The group H has three
generators
$$
a=((1,0),0); \quad b=((0,1),0); \quad c=((0,0),1)
$$
and is defined by the relations
$ [a,b]=e; \quad [b,c]=e; \quad
[c,a]=b, $ where $[\cdot,\cdot]$ stands for the commutator, and $e=((0,0),0)$ is the unit of the group.

{\bf Lemma 1.} {\it For each $N \in \mathbb N$ there exists an automorphism
 $\phi$ of the group $H$ such that $R(\phi)=2N$.}

{\bf Proof} Let $N \in \mathbb N$ be a given integer. Consider the following automorphism:
$$
\phi((m,k),s)=\Big ( \big
(m(1-N)+sN,-k+\frac{m(m-1)}{2}(1-N)-\frac{s(s-1)}{2}N+msN
\big),m-s\Big ).
$$
To check that this mapping defines an automorphism, one can either verify directly that it is a bijective
endomorphism or use the description of all Heisenberg group automorphisms given in [8].

Denote
$$
\quad \quad \quad \quad
Q_0(m,s)=\frac{m(m-1)}{2}(1-N)+\frac{s(s-1)}{2}N+msN.
$$
Then
$$
\phi((m,k),s)=\big ((m(1-N)+sN,-k+Q_0(m,s)),m-s\big ).
$$

Describe the $\phi$- conjugacy class of an arbitrary element
$h=((m,k),s)$ of the group $H$. Let $g=((g_1,g_2),g_3)$. Then
$g^{-1}=((-g_1,g_1g_3-g_2),-g_3)$ and
$$
\phi(g^{-1})=\big( (-g_1(1-N)-g_3N,g_2-g_1g_3+Q_0(-g_1,-g_3)), g_3-g_1 \big).
$$
Assuming $-g_1g_3+Q_0(-g_1,-g_3)=:Q_1(g_1,g_3)$, we have
\begin{multline*}
g h \phi(g^{-1})=((g_1,g_2),g_3)\cdot ((m,k),s) \cdot ( (-g_1(1-N)-g_3N,g_2+Q_1(g_1,g_3)), g_3-g_1)=\\
=((g_1+m,g_2+k+g_3m),g_3+s)\cdot( (-g_1(1-N)-g_3N,g_2+Q_1(g_1,g_3)), g_3-g_1)=\\
=((g_1+m-g_1+g_1N-g_3N,2g_2+k+Q_2(g_1,g_3,s,m)),s+2g_3-g_1)=\\
=((m+(g_1-g_3)N,k+2g_2+Q_2(g_1,g_3,s,m)),s+2g_3-g_1),
\end{multline*}
where $Q_2(g_1,g_3,s,m)=g_3m+Q_1(g_1,g_3)-g_3g_1(1-N)-sg_1(1-N)-g_3^2N-sg_3N$.

Thus, the $\phi$-conjugacy class of the element $h=((m,k),s)$ has
the form
$$
\begin{aligned}
{h}_{\phi}=&\{g h \phi(g^{-1})\mid h \in H\}=\\
=&\{\big((m+(g_1-g_3)N,k+2g_2+Q_2(g_1,g_3,s,m)),s+2g_3-g_1\big)\mid g_1,g_2,g_3\in \mathbb Z \}.
\end{aligned}
$$
This class is parameterized by a triplet $g_1,g_2,g_3$ of integers.
Perform the following re\-para\-me\-ter\-i\-za\-tion:
$$
 g_1-g_3=:f_1;\quad g_2=:f_2; \quad 2g_3-g_1=:f_3.
$$
This reparametrization is valid because its determinant equals $-1$,
and we have
$$
\{((m,k),s)\}_{\phi}=\{ \big( (m+f_1N,k+2f_2+Q_3(f_1,f_3,s,m)),s+f_3 \big) \mid f_1,f_2,f_3\in \mathbb Z  \}.
$$

The latter implies that if $m$-components of two group elements do
not coincide modulo $N$, then these elements are not equivalent
(i.e., not $\phi$-conjugate). Assume
$$
H_r=\{ \big( (m,k),s\big) \mid m \equiv r (\bmod N) \}; \quad
r=\overline{0,\ldots,N-1}.
$$
Then two elements from distinct sets $H_r$ can never be equivalent to each other. It is also evident that these
sets form a partition of $H$.

Fix $r$ and consider equivalence classes of the elements $((r,0),0)$ and
$((r,1),0)$.
\begin{multline}
\label{last}
\begin{aligned}
&\{((r,0),0)\}_{\phi}=\{ \big( (r+f_1N,2f_2+Q_3(f_1,f_3,r,0)),f_3 \big) \mid f_1,f_2,f_3\in \mathbb Z  \};\\
&\{((r,1),0)\}_{\phi}=\{ \big( (r+f_1N,2f_2+1+Q_3(f_1,f_3,r,0)),f_3 \big) \mid f_1,f_2,f_3\in \mathbb Z  \}.
\end{aligned}
\end{multline}
Now take an arbitrary element $((m_0,k_0),s_0),\quad m_0
\equiv r (\bmod N)$ of the set $H_r$ and show that it belongs
to one of these two classes.

Indeed, one can choose $f_1$ and $f_3$ in (\ref{last}) so that
$r+f_1N=m_0, f_3=s_0$. Thereby, the value
$Q_3(f_1,f_3,r,0)$ is fixed. Depending on the parity of this value, one can choose element
 $f_2$ so that either $2f_2+Q_3(f_1,f_3,r,0)=k_0$,
or $2f_2+Q_3(f_1,f_3,r,0)=k_0+1$. In the first case $((m_0,k_0),s_0)
\sim ((r,0),0)$ and in the second case $((m_0,k_0),s_0) \sim
((r,1),0)$.

Thus, the group $H$ is partitioned into $N$ subsets $H_r$ such that
elements from different subsets are not equivalent and each $H_r$
is partitioned into exactly two equivalence classes. Therefore,
$R(\phi)=2N$ and we are done.\\

\begin{center}{\bf 2. DETERMINATION OF FIXED REPRESENTATIONS}\\   \end{center} Describe
finite-dimensional irreducible unitary representations
 of the discrete Hei\-sen\-berg group
$H=\mathbb{Z}^2 \rtimes \mathbb{Z}$, using the technique of induced representations.
Afterwards, we find fixed representations for a specific automorphism with the Reidemeister number 2.

The dual object for $\mathbb Z ^2$ is the torus $\T^2=\mathbb R^2 /
\mathbb Z^2$. A pair  $\chi=(\xi, \eta) \in \T^2$ corresponds to the
character  $(m,k) \mapsto e^{2\pi i (m\xi+k\eta)}$. This torus is a
right $H$-space with respect to the action $ \chi h
(m,k)=\chi(h*((m,k),0)*h^{-1}). $ The action for $((m,k),s)$ is
defined by the formula $
(\xi, \eta) \mapsto (\xi,\eta) \bigl(\begin{smallmatrix} 1 & 0\\
s& 1 \end{smallmatrix} \bigr)=(\xi+s\eta, \eta). $

In accordance with the Mackey theorem [10], each irreducible
representation $\rho$ of the group $H$ has the form $\Ind(H, Y,
\beta)$, where $Y$ is the stabilizer of some point $\chi \in
\widehat N$, and the restriction of $\beta$
 onto $N$ is scalar and, moreover, is a multiple of the character $\chi$.
Conversely, due to [11, p. 159, Theorem 5], if $\beta$
is
an irreducible representation of $Y$, whose restriction onto $N$ is a multiple of a character, then
 $\Ind(H, Y, \beta)$ is irreducible. Finally, the theorem from [12, p. 155] implies that each
orbit $O$ is homeomorphic to the quotient
of $H$ over the stabilizer of any point from the orbit.

Describe the finite-dimensional irreducible unitary representations
of the discrete Hei\-sen\-berg group $H$. Due to aforesaid, each such a
representation is induced from the stabilizer of some point
$\chi=(\xi, \eta)\in \T^2$, whose orbit $\widehat O_\chi$ is finite.
Let the orbit of some point $\chi$ consist of $p$ points. Then the
stabilizer of the point $\chi$ is the subgroup
$$
Y=\{((m,k),ps) \mid m,k,s \in \mathbb Z \} = \mathbb Z ^2 \rtimes
(p\mathbb Z).
$$

In accordance with [11, p. 158, Lemma 3], the multiplication by the character
cor\-res\-pon\-d\-ing to $\chi$, is a
bijection between the sets of irreducible unitary representations of the groups $p\mathbb Z$ and $Y=\mathbb Z ^2 \rtimes
(p\mathbb Z)$.

Therefore, any  $p$-dimensional irreducible unitary representation
of the group $H$ ($p<\infty$) can be obtained by means of the
following algorithm: 1) choose a point $\chi=(\xi, \eta) \in \T^2$,
having the orbit of cardinality $p$; 2) choose an irreducible
representation of the subgroup $p\mathbb Z = \{((0,0),ps)\mid s \in
\mathbb Z\}$; 3) multiply this representation by the character
$\chi$ and obtain a representation $\beta$ of the corresponding
subgroup $Y$; 4) form the representation $\rho$ of the group $H$,
induced by $\beta$.

Choose an arbitrary number $\alpha \in [0,1)$ and hence an arbitrary irreducible representation
 $\pi$ of the subgroup $(p\mathbb Z)$:
$$
\pi((0,0),ps)=e^{2\pi i s \alpha}.
$$

Multiplying by the character $\chi=(\xi, \eta)$ we obtain a representation $\beta$
of the subgroup $Y$:
$$
\beta((m,k),ps)=\chi(m,k)e^{2\pi i s \alpha}=e^{2\pi
i(m\xi+k\eta+s\alpha)}.
$$
Now we have to form the representation $\rho$ of the group $H$,
induced by the representation $\beta$. To do that, produce a
realization of an induced representation in the space of
$L^2$-functions on the homogeneous space
 $X=Y \diagdown H$ (see [10]) in the case of a discrete group.

Let $H$ be a discrete group, $Y$ be its subgroup, $\beta$ be a unitary representation of $Y$ in a Hilbert space
 $V$, $X=Y \diagdown H$ be the corresponding right homogeneous space.
Fix a mapping $s \colon X \to G$, possessing the
property $s(Hg) \in Hg$. Then the induced representation $\rho$ is given in the space $L^2(X,V)$ by the formula \\
$ [\rho(h)f](x)=A(h,x)f(xh), $ where the operator-valued function
$A(h,x)$ is given by the equality $ A(h,x)=\beta(y), $ where the
element $y \in Y$ is defined by the relation $ s(x)h=ys(xh). $

In our case the space of the representation $\beta$ is one-dimensional and hence
 $A(h,x)$ is a complex-valued function. Note that
$Y((m,k),s)=Y((0,0),s\bmod p)$, and, therefore,
$$
Y \diagdown H=\{Y((0,0),0),\quad Y((0,0),1), \quad\ldots ,\quad Y((0,0),p-1) \}=:\{x_0,x_1,\ldots, x_{p-1}\}.
$$
Choose the mapping $s\ \colon X \to G$ by assuming $ s \colon Y((m,k),s)
\longmapsto ((0,0), s\bmod p). $

Now we have to determine an element $y \in Y$ using the relation $s(x)h=ys(xh)$   and starting from given
elements $x \in X=Y \diagdown H$ and $h \in H$. Let
$h=((m,k),s)$, $x=Y((0,0),j).$ Then $s(x)=((0,0),j)$ and
$$
 s(x)h=((0,0),j)((m,k),s)=((m,k+jm),s+j);
$$
$$
s(xh)=s(Y((0,0),j)((m,k),s))=s(Y(m,k+jm),s+j)=((0,0),(s+j)\bmod p).
$$
Suppose the element $y$ looked for has the form $y=((y_k,y_m),y_s).$ In this case we have
$$
ys(xh)=((y_k,y_m),y_s)((0,0),(s+j)\bmod p)=((y_m,y_k),y_s+(s+j)\bmod p).
$$
Therefore,
$$
y_m=m, \quad y_k=k+jm, \quad y_s=s+j-(s+j)\bmod p.
$$
Assuming $[l]_p:=l-l\bmod p$, we have $y=((m,k+jm),[s+j]_p)$ and
$$
A(h,x)=\beta(y)=\beta((m,k+jm),[s+j]_p)=e^{2\pi i
(m\xi+(k+jm)\eta+\frac{[s+j]_p}{p})}=
$$
$$
=e^{2\pi i(m\xi+(k+jm)\eta+[\frac{s+j}{p}])},
$$
where $[r]$ denotes the integer part of $r$.

At last, the induced representation of $\rho$ in the space
$L^2(X)=L^2(\{x_0,x_1,\ldots, x_{p-1}\})$ is given by the formula
$$
[\rho(h)f](x)=A(h,x)f(xh)=e^{2\pi i(m\xi+(k+jm)\eta+[\frac{s+j}{p}])}f(xh);
$$
$$
[\rho((m,k),s)f](x_j)=e^{2\pi i(m\xi+(k+jm)\eta+[\frac{s+j}{p}])}f(x_{(j+s)\bmod p});
\quad j=\overline{0,p-1}.
$$
In the space $L^2(X)=L^2(\{x_0,x_1,\ldots, x_{p-1}\})
\cong \mathbb C ^p$ choose an orthonormal basis $\epsilon_0, \epsilon_1,
\ldots,$ $\epsilon_{p-1}$, where $\epsilon_j$
is the indicator function of the point $x_j \in X$. In this basis our representation is given by the formula
\begin{equation}
\label{rho}
\rho((m,k),s)\colon \epsilon_{j}\mapsto
e^{2\pi i(m\xi+(k+jm)\eta+[\frac{s+j}{p}]\alpha)}\epsilon_{(j-s)\bmod p};
\quad j=\overline{0,p-1}.
\end{equation}
Thus, all $p$-dimensional irreducible unitary representations of the group $H$ have the form (\ref{rho})
for some numbers $\xi, \eta, \alpha \in [0,1)$, such that the orbit  $\chi=(\xi,\eta) \in \T^2$
consists of $p$ points. The action of  $((m,k),s)$ on $\T^2$ is given
by the formula $ (\xi, \eta) \mapsto (\xi+s\eta, \eta).$ Therefore, the orbit of this action has cardinality $p$
if and only if $\eta$ is an irreducible fraction with the denominator $p$. Thus, the following statement is proved.

{\bf Lemma 2.} {\it All $p$-dimensional irreducible unitary representations of the discrete
Hei\-sen\-berg group $H$
have the form(\ref{rho}), where $\epsilon_j$ are the elements of an orthonormal basis of the
 $p$-dimensional space, $\xi, \eta, \alpha \in [0,1)$,
and $\eta$ is an irreducible fraction of denominator $p$.}

Find the character $\chi_{\rho}$ of the representation (\ref{rho}). The matrix of the operator $\rho((m,k),s)$
is diagonal for $s\equiv 0 \bmod p$ and for other $s$ its diagonal entries are zero. Therefore, if $s$ is not divisible by
$p$, then we have $\chi_{\rho}((m,k),s)=0$. Thus,
$$
\chi_{\rho}((m,k),s)=
\delta_{s\bmod p}^0\sum_{j=0}^{p-1}
\exp(2\pi i(m\xi+k\eta+jm\eta+[\frac{s+j}{p}]\alpha)).
$$
For $s$ divisible by $p$ and $j \in \overline{0,p-1}$ the relation
$[\frac{s+j}{p}]=\frac{s}{p}$ is valid. Further, taking the
multiples which do not contain $j$ outside of the summation, we get
$$
\chi_{\rho}((m,k),s)=\delta_{s\bmod p}^0 \exp  \big ( 2\pi
i(m\xi+k\eta+\frac{s}{p}\alpha)\big ) \sum_{j=0}^{p-1}e^{2\pi i m
\eta j}.
$$
The resulting sum is the sum of the first $p$ terms of a geometric series. Note that
 $e^{2\pi i m \eta}=1$ if and only if
$m \eta \in \mathbb Z$, i.e., iff $m\equiv 0 \bmod p$.
Therefore,
$$
\sum_{j=0}^{p-1}e^{2\pi i m \eta j}=\begin{cases}
p, &\mbox{ if }m \equiv 0\bmod p;\\
\frac{\exp(2\pi i m\eta p)-1}{\exp(2\pi i m\eta )-1}
&\mbox{ if }m \not \equiv 0\bmod p.
\end{cases}
$$
But $\eta p \in \mathbb Z$. Hence $\exp(2\pi i m\eta p)-1=0$. As a result, we get
$$
\chi_{\rho}((m,k),s)=\begin{cases} p\cdot e^{2\pi i
(m\xi+k\eta+\frac{s}{p}\alpha)}, &\mbox{ if }s \equiv 0\bmod p
\mbox{~ and } m \equiv 0\bmod p;\\
0 &\mbox{ otherwise}.
\end{cases}
$$

Consider the following automorphism  $\phi$ of the group $H$:
$$
\phi((m,k),s)=\Big ( \big (s+m,-k+\frac{m(m-1)}{2}+sm \big ),m \Big
).
$$
By reasoning similar to those used in the proof of Lemma 1, it is not difficult to see that $R(\phi)=2$.
Calculate the character of the representation $\rho\phi$:
\begin{multline*}
\rho\phi((m,k),s)\epsilon_j=\rho((s+m,-k+\frac{m(m-1)}{2}+sm),m)\epsilon_j=\\
=\exp2\pi i ((s+m)\xi+(-k+\frac{m(m-1)}{2}+sm+js+jm)\eta+[\frac{m+j}{p}]\alpha)
\epsilon_{(j-m)\bmod p}.
\end{multline*}
As in the previous calculations, the character $\chi_{\rho\phi}$
vanishes on the elements which have $m \not \equiv 0 \bmod p$. Then
we use the fact that $[\frac{m+j}{p}]=\frac{m}{p}$,
if $m\equiv 0 \bmod p$ and $j\in \overline{1,p}$. We have:
$$
\chi_{\rho\phi}((m,k),s)= \delta_{m\bmod p}^0\exp2\pi i ((s+m)\xi
+(-k+\frac{m(m-1)}{2}+sm)\eta+\frac{m}{p}\alpha)\sum_{j=0}^{p-1}
e^{2\pi i (s+m)\eta j}
$$
and
$$
\sum_{j=0}^{p-1}e^{2\pi i (s+m)\eta j}=\begin{cases}
 p, &\mbox{ if} (s+m)\eta \in \mathbb Z \mbox{,~~i.e. if }(s+m)\equiv 0 \bmod p,\\
 0, &\mbox{ if } (s+m)\eta \not\in \mathbb Z.
\end{cases}
$$
So,
$$
\chi_{\rho\phi}((m,k),s)=\begin{cases}
p\cdot e^{2\pi i ((s+m)\xi+(-k+\frac{m(m-1)}{2}+sm)\eta+\frac{m}{p}\alpha)}, &
\mbox{ if }s \mbox{ and } m \equiv 0\bmod p,\\
0, &\mbox{ otherwise}.
\end{cases}
$$

Let us find finite-dimensional fixed points of the mapping $\wh\phi$. The trivial one-dimensional representation is one
of them. To find another one, write down the condition for characters coinciding for the representations $\phi$
and $\rho\phi$:
$$
e^{2\pi i
((s+m)\xi+(-k+\frac{m(m-1)}{2}+sm)\eta+\frac{m}{p}\alpha)}=
e^{2\pi i (m\xi+k\eta+\frac{s}{p}\alpha)} \mbox{ для }s, m \equiv
0\bmod p.
$$
Assume $p=2$, i.e., we look for a fixed representation among two-dimensional ones. Тогда $\eta=\frac{1}{2}$ and
$e^{2\pi i ((s+m)\xi+\frac{1}{2}(-k+\frac{m(m-1)}{2}+sm)+
\frac{m}{2}\alpha)}=e^{2\pi i
(m\xi+\frac{1}{2}k+\frac{s}{2}\alpha)} \mbox{ for even }s, m.
$

Assuming $s=:2t, m=:2q$,  we obtain:
$
e^{2\pi i ((2t+2q)\xi+\frac{1}{2}(-k+q(2q-1)+4tq)+q\alpha)}=
e^{2\pi i (2q\xi+\frac{1}{2}k+t\alpha)}$   for any $t, q \in
\mathbb Z.
$
This condition is equivalent to
$$
\big ((2t+2q)\xi+\frac{1}{2}(-k+q(2q-1)+4tq)+q\alpha\big)- \big
(2q\xi+\frac{1}{2}k+t\alpha\big)\in \mathbb Z \mbox{ для любых }t,
q \in \mathbb Z.
$$
Collecting terms and throwing out the terms $q^2, 2tq, -k$, which
are certainly integer, we get the following equivalent condition:
$$
\begin{aligned}
&2t\xi-\frac{q}{2}+q\alpha -t\alpha \in \mathbb Z \mbox{ для любых
}t, q \in \mathbb Z, or\\
&t(2\xi-\alpha)+q(\alpha-\frac{1}{2})\in \mathbb Z, \forall t, q \in
\mathbb Z.
\end{aligned}
$$
Evidently, this condition holds for $\alpha=\frac{1}{2},
\xi=\frac{1}{4}$.

Thus, the second fixed representation is two-dimensional and it is given by formula (\ref{rho}) with $\alpha=\frac{1}{2},
\xi=\frac{1}{4}, \eta = \frac{1}{2}$, i.e.
$$
\rho_2((m,k),s)\colon \epsilon_{j}\mapsto
e^{2\pi i(\frac{m}{4}+\frac{k+jm}{2}+
\frac{1}{2}[\frac{s+j}{2}])}\epsilon_{(j-s)\bmod 2};\quad j=0,1.
$$
where $\epsilon_1, \epsilon_2$ -is an orthonormal basis of the two-dimensional space.\\

\begin{center} {\bf 3. WREATH PRODUCTS OF ABELIAN GROUPS AND $\mathbb
Z$}\\   \end{center}
Construct a class of groups which are finitely generated, possess the RP property, but are not almost polycyclic.
Namely, we prove that for an arbitrary finitely generated Abelian group $A$, the wreath product
$A \wr\mathbb Z$ possesses the RP property.

To do this, we show that some subgroup $\mathcal A \subset A
\wr\mathbb Z$ is {\it characteristic} (i.e., it is invariant under all
automorphisms of the group$A\wr \mathbb Z$),
and then apply Theorem 3.10 from [5]. At last, we show that any such
wreath product is not an almost polycyclic group.

Let two groups $A$ and $B$ be given. The {\em wreath product} of these groups is defined
as the following semidirect product:
$$
A\wr B := \left( \bigoplus_{i\in B}A_i \right ) \rtimes B
$$
where each subgroup $A_i$ is isomorphic to $A$ and the action of $B$ on
$\bigoplus A_i$ is given by the formula $b \cdot A_i = A_{bi} $.

Recall that, as for any semidirect product, the action of an element $b \in B$ on the semidirect multiple
$\bigoplus A_i$ coincides with the inner automorphism of the product $A \wr B =\left(
\bigoplus A_i \right ) \rtimes B$, given by $b$.

In what follows, by $A$ we always denote a finitely generated Abelian group and $B$
coincides with the group of integers. By $\mathcal A$
we denote the sum $\bigoplus_{i\in \mathbb Z}A_i$ from the definition of a wreath product.

It is not difficult to see that the group $G= A \wr \mathbb Z$ is generated by generators of the
groups $A$ and $\mathbb Z$. Therefore, $G$ is finitely generated.

Derive the statement concerning the characteristic property of the subgroup
$A$ in the particular case when $A$ is a free Abelian group. Thus,
let $A=\mathbb Z ^k$. We naturally identify the Abelian group
$\mathcal A = \left(\bigoplus_{i\in \mathbb Z}A_i \right )$ with the group
$\mathbb Z [x,x^{-1}]$ of Laurent polynomials with respect to addition. Namely, a given sequence
$$\ldots
0,(a^{(N)}_0,\ldots a^{(N)}_{k-1}),\cdots (a^{(M)}_0,\ldots
a^{(M)}_{k-1}),0 \ldots
$$
of elements from $\mathbb Z^k$ can be considered (omitting the parenthesis) as a sequence of integers identified with
coefficients of the corresponding Laurent polynomial $P(x)$:
$$
P(x)=a^{(N)}_0 x^{Nk}+a^{(N)}_1 x^{Nk+1}+ \cdots
+a^{(M)}_{k-1}x^{Mk+(k-1)}.
$$

For example, if $A=\mathbb Z^2$, and $\bar a = (a_i)_{i \in \mathbb Z}
\in \mathcal A$ has $a_{-1}=(3,-1);\quad a_{0}=(5,7)$, and all others
$a_i$-s are zero, then $\bar a$ corresponds to the polynomial $3x^{-2}-x^{-1}+5+7x. $

Thus, an arbitrary element of the group $G=\mathbb Z^k \wr
\mathbb Z = \mathcal A \rtimes \mathbb Z$ may be written down in the form
$(P(x),n)$, where $P(x) \in \mathcal A$ is a Laurent polynomial and, а $n$
is an integer.

From the definition of the wreath product, it is clear that the action of the unit (the generator of
$\mathbb Z$) on $\bigoplus_{i\in \mathbb Z}A_i$ is the multiplication of the polynomial by $x^k$.

Therefore, the group multiplication is given by the formula
$$
(P(x),r)(Q(x),s)=(P(x)+x^{rk}Q(x),{r+s}),
$$
and the inverse element is given by the formula $(P(x),r)^{-1}=(-x^{-kr}P(x),{-r}).$

Consider the following set of mappings $\epsilon_j \colon \mathcal A
\to \mathbb Z;\quad j=\overline{0,k-1}$;
$$
\epsilon_j \colon \sum a_nx^n \mapsto \sum_{n\equiv j \pod {k}}a_n.
$$

{\bf Lemma 3.} {\it The quotient group of $G=\mathbb Z^k \wr \mathbb
Z$ over its commutant is isomorphic to $\mathbb Z^{k+1}$, and the
canonical projection is given by the formula $
p(P(x),r)=(\epsilon_0(P(x)),\ldots,\epsilon_{k-1}(P(x)),r). $ } {\bf
Proof.} First, notice that the commutant of $G$ lies inside of
$\mathcal A = \bigoplus_{i\in \mathbb Z}A_i $, as the quotient
$G/\mathcal A=\mathbb Z$ is Abelian.

Calculate explicitly the commutator of two arbitrary elements of the group
$$
\left [ (P(x),r),(Q(x),s) \right ]
=((1-x^{ks})P(x)+(x^{kr}-1)Q(x),0).
$$

It is not difficult to see that $\epsilon_j\big
((1-x^{ks})P(x)+(x^{kr}-1)Q(x)\big )=0, $ for any $j$ because the multiplication of
the polynomial by a $k$-multiple degree of $x$  does not change the values of
$\epsilon_j$. Hence, the commutant contains
only polynomials such that $\epsilon_j$ are equal to zero on them. Moreover, the commutant contains all such
polynomials. Indeed, taking $P(x)=x^n, r=0, Q(x)=0$, we easily obtain that $G'$ contains all polynomials of
the form $x^n-x^{n+k}$, and any polynomial with zero $\epsilon_j$ - s can be represented as an integer combination of the
latter polynomials.

Thus,
$
G'=\{ (P(x),0) \mid P(x) \in \mathbb Z [x], \epsilon_0(P)=0,
\ldots,\epsilon_{k-1}(P)=0\}
$
which easily implies the assertion of the Lemma.

{\bf Lemma 4.} {\it The subgroup $\mathcal A = \bigoplus_{i \in \mathbb
Z} (\mathbb Z ^k)$ is characteristic in $\mathbb Z^k \wr \mathbb
Z$.}

{\bf Proof}
Consider the following generators system of the group $\mathbb Z^k \wr
\mathbb Z$:
$$
a_0=(1,0),\quad a_1=(x,0),\quad \ldots\quad, \quad
a_{k-1}=(x^{k-1},0)\mbox{ --- the generators of } \mathbb Z^k;\quad
b=(0,1).
$$
Let $\phi$ be an automorphism of $G$ and $ \phi(a_i)=(P_i(x),{s_i}),
i=\overline{0,k-1} ; \quad \phi(b)=(Q(x),r).$ It is sufficient
to show that $s_i=0, i=\overline{0,k-1}$. Indeed, each element of $\mathcal A$ is a finite product of elements
of the form $\alpha=b^t \cdot a \cdot b^{-t},
 a=a_0^{u_0}\cdots a_{k-1}^{u_{k-1}}$, and if all $s_i=0$, then
$
\phi(\alpha)=(\phi(b))^t \phi (a) \phi
(b)^{-t}=(*,rt)(*,0)(*,-rt)=(*,0) \in \mathcal A.
$

Note that the automorphism $\phi$ induces the automorphism of the quotient $G/G'
\cong \mathbb Z^k\oplus \mathbb Z$, having the
following matrix:
$$\Pi=
\begin{pmatrix}
\epsilon_0(P_0) & \epsilon_0(P_1) & \ldots & \epsilon_0(P_{k-1}) &
\epsilon_0(Q)\\ \epsilon_1(P_0) & \epsilon_1(P_1) & \ldots &
\epsilon_1(P_{k-1}) & \epsilon_1(Q)\\ \vdots & \vdots & \ddots &
\vdots & \vdots\\ \epsilon_{k-1}(P_0) & \epsilon_{k-1}(P_1) & \ldots
& \epsilon_{k-1}(P_{k-1}) & \epsilon_{k-1}(Q)\\ s_0 & s_1 & \ldots &
s_{k-1} & r
\end{pmatrix}
$$

The generators $a_i$ lie in an Abelian subgroup and hence commute. Therefore, their
$\phi$-images also commute and hence
$\phi(a_i)\phi(a_j)=\phi(a_j)\phi(a_i)$,
which implies
$$
\begin{aligned}
&(P_i(x)+x^{k s_j}P_j(x),{s_i+s_j})=(P_j(x)+x^{ks_i}P_i(x),{s_i+s_j});\\
&P_i(x)(1-x^{k s_j})=P_j(x)(1-x^{k s_i}); (3)
\end{aligned}
$$
So if $s_j=0$ then either $s_i=0$ or $P_j(x)=0$. But the latter is impossible because the element
 $\phi(a_j)=(P_j(x),{s_j})$ cannot be equal to the group unit.
Hence, $s_j=0$ implies $s_i=0$. Due to the arbitrariness of $i$, $j$ the latter implication means
that either all numbers $s_i$ are equal to zero or none of them is equal to zero. In the first
case the proof is completed. Prove that the second case cannot be valid.

Assume that $s_i \ne 0$, for all $i$. In this case we can divide the polynomial equality $(3)$ by $(1-x^k)$.
Dividing and applying $\epsilon_{k-1}$ to both parts of the equality, we get
$$
\begin{aligned}
&P_i(x)(1+x^k+\ldots+x^{k(s_j-1)})=P_j(x)(1+x^k+\ldots+x^{k(s_i-1)});\\
&\epsilon_{k-1}(P_i)\cdot s_j = \epsilon_{k-1}(P_j)\cdot s_i;
\quad (4)
\end{aligned}
$$

Similar to the above reasoning $a_i$ and $ba_i b^{-1}$
lie in an Abelian subgroup, hence their images commute.
We have
$$
\phi(ba_i b^{-1})=(Q(x),r) \cdot (P_i(x),{s_i}) \cdot (-x^{-kr}Q(x),{-r})
=((1-x^{ks_i})Q(x)+x^{kr}P_i(x),{s_i}).
$$
The commutation condition $\phi(a_i)\cdot \phi(ba_i b^{-1}) =\phi(ba_i
b^{-1})\phi(a_i)$ may be evidently transformed to the following form:
$$
P_i(x)(1-x^{ks_i})(1-x^{kr})=Q(x)(1-x^{ks_i})^2.
$$
Since we assume that $s_i \ne 0$ for any $i$, then we can divide the latter equation by  $(1-x^{ks_i})(1-x^k)$. We get
$$
P_i(x)(1+x^k+\ldots+x^{k(r-1)})=Q(x)(1+x^k+\ldots+x^{k(s_i-1)});
$$
Apply $\epsilon_{k-1}$ to this equation:
$$
\epsilon_{k-1}(P_i)\cdot r = \epsilon_{k-1}(Q)\cdot s_i. \quad (5)
$$
Combining $(4)$ and $(5)$, we see that all lower $(2\times 2)$-minors of the matrix $\Pi$ are degenerate and so the matrix
is degenerate. But it is impossible because $\Pi$ is the matrix of an automorphism of the group $\mathbb Z ^ {k+1}$.
The contradiction obtained proves that all $s_i$ are zero, and we are done.

Now let us proceed from the case $A=\mathbb Z ^ k$ to a more general case.

{\bf Lemma 5.} {\it Let $A$ be a finitely generated Abelian group. Then the subgroup $\mathcal A = \bigoplus A_i$
is characteristic in the group $A \wr \mathbb Z = \left( \bigoplus A_i
\right ) \rtimes \mathbb Z$.}

{\bf Proof}
We have $A_i=(\mathbb Z ^k)_i \oplus T_i$, where $T_i$ is the torsion subgroup of $A_i$.
Write down elements of $\mathcal A$ as
pairs $(\bar b,\bar f)$, where $\bar b \in \bigoplus( \mathbb Z^k)_i$, and
$\bar f \in \bigoplus T_i$.

Show that the subgroup $\mathcal T= \bigoplus T_i$ is characteristic
in $A \wr \mathbb Z$. First we suppose that there exists $g \in
\mathcal T :\quad \phi(g)= \phi((\bar 0, \bar f),0)=((*,*),d)$,
where $d\ne 0$. Let $N$ be the maximal order of elements from $T$.
Then, on the one hand, $\phi(g^N)=\phi(e)=e$, and, on the other
hand, $\phi(g^N)=((*,*),{Nd})\ne e$.

So, $d=0$ and $\phi((\bar 0, \bar
f),0)=((\bar{b_1},\bar{f_1}),0)$. Now the same reasoning proves that $\bar{b_1}=\bar 0$, as
$((\bar{b_1},\bar{f_1}),0)^N$=$((N\bar{b_1},N\bar{f_1}),0)$.

Thus, the subgroup $\mathcal T$ is characteristic and for any $\phi \in \Aut(A \wr \mathbb Z)$ we have the following commutative
diagram:

$$
\begin{CD}
0 @ >>> \mathcal T @>i>> A \wr \mathbb Z @>p>> \mathbb Z^k \wr
\mathbb Z @>>>
0\\
@. @VV\phi' V @VV\phi V @ VV\overline{\phi}V  @.\\
0 @ >>> \mathcal T @>i>> A \wr \mathbb Z @>p>> \mathbb Z^k \wr
\mathbb Z@>>>
0\\
\end{CD}
$$

Assume that the statement of the lemma is not valid. Then there
exists an element $g=((b,f),0)$ of the subgroup $\mathcal A$, such
that $\phi(g)\not \in \mathcal A$. We have
$\phi(g)=\phi((b,0),0)\cdot\phi((0,f),0) \not \in \mathcal A$;
hence, $\phi((b,0),0) \not \in \mathcal A$, because $\phi((0,f),0)
\in \mathcal T  \subset \mathcal A$ according to what was proved.

The fact that $\phi((b,0),0) \not \in \mathcal A$, implies $p \circ
\phi((b,0),0)\notin \oplus(\mathbb Z^k)_i$, which is equivalent to
$\overline{\phi} \circ p((b,0),0)\notin \oplus(\mathbb Z^k)_i$.
Тhus, the subgroup $\bigoplus (\mathbb Z^k)_i$ is not invariant with
respect to the automorphism $\bar {\phi}$ of the group группы
$\mathbb Z^k \wr \mathbb Z$. But this contradicts the previous
lemma, and we are done.

Theorem 3.10 from [5] states that if in a group extension
$$0\to H \to G \to G/H \to 0$$ the subgroup $H$ is characteristic
and has the  RP property, and $G/H$ is a finitely generated group
with finite conjugacy classes, then $G$ possesses the RP property.

{\bf Theorem 1.} {\it For an arbitrary finitely generated Abelian
group $A$, the wreath product $A \wr\mathbb Z$ possesses the RP
property.}

{\bf Proof.} Show that the group extension $0\to \mathcal A \to A
\wr \mathbb Z \to \mathbb Z \to 0$ satisfies the conditions of the
theorem presented above.

The subgroup $\mathcal A$ is Abelian and, therefore,
it possesses the RP-property ([5]). Further, we have proved that this
subgroup is characteristic in the group $A\wr \mathbb Z$. At last,
the quotient $(A\wr\mathbb Z)/\mathcal A = \mathbb Z$ has finite
conjugacy classes (consisting of one element). The theorem is proved.

The latter result is of particular interest because such groups gave
the first known example of finitely generated RP groups being not
almost polycyclic.

Recall that a group $G$ is said to be \emph{polycyclic}, if there
exists a sequence of subgroups (\emph{a polycyclic series}) $
\{e\}=G_0  \lhd G_1 \lhd \cdots \lhd G_k=G, $ such that all
$G_{i+1}/G_i$ are cyclic. A group is said to be \emph{almost
polycyclic}, if it contains a polycyclic subgroup of finite index.

Show that for an arbitrary finitely generated Abelian group $A$, the
group $A \wr \mathbb Z$ s not almost polycyclic. Indeed, the main
theorem from [13] immediately implies that each characteristic
subgroup of an almost polycyclic group is finitely generated. But we
have proved that the subgroup $\mathcal A = \bigoplus_{i\in \mathbb
Z}A_i$ of the group $A \wr \mathbb Z$ -is characteristic and the
fact that it is infinitely generated is evident.\\ \\ {\it The work
was partly supported by the Russian Foundations for Basic Research
(RFFI), project 05–01–00923-a, and by the Program “Universities of
Russia.”''}.\\

\begin{center} {\sc BIBLIOGRAPHY} \end{center}

\small {
\begin{enumerate}

\vspace{-0.2cm}
\item{ \emph {A.~L. Fel'shtyn}, {The {R}eidemeister number of any automorphism of a
  {G}romov hyperbolic group is infinite} // Zap. Nauchn. Sem. S.-Peterburg.
  Otdel. Mat. Inst. Steklov. (POMI) \textbf{279} (2001), no.~6 (Geom. i
  Topol.), 229--240, 250. }\\

\vspace{-0.65cm}
\item{\emph{Gon{\c{c}}alves D., Wong P.} {Twisted conjugacy classes in
exponential growth groups} // Bull. London Math. Soc. 2003.\textbf{35}, N.~2 261--268.}\\ 

\vspace{-0.65cm}
\item{\emph{{Fel'shtyn A., Gon{\c{c}}alves D.}} Reidemeister numbers
of Automorphisms of {B}aumslag-{S}olitar groups // Algebra and
discrete Mathematics 2006.\textbf{3} 36--48.}\\

\vspace{-0.65cm}
\item{\emph{Gon{\c{c}}alves D., Wong P.} {Twisted conjugacy classes in
wreath products} // Internat. J. Algebra and Comput. 2006.\textbf{16}, N 2 875--886.}\\

\vspace{-0.65cm}
\item{\emph{Fel'shtyn A., Troitsky E.}. {Twisted {B}urnside theorem} //
Preprint MPIM2005-46,
  Max-Planck-Institut f{\"u}r Mathematik, 2005.}\\ 

\vspace{-0.65cm}
\item{\emph{Serre J.-P.} {Linear Representations of Finite Groups}. Springer Verlag, Berlin Heidelberg- New York,
1977}\\

\vspace{-0.65cm}
\item{\emph{Fel'shtyn A., Troitsky E.}. {A twisted {B}urnside theorem for
countable groups and {R}eidemeister numbers} // Proc. Workshop
Noncommutative Geometry and Number
  Theory (Bonn, 2003) (Ed. by K.~Consani, M.~Marcolli, and Yu. Manin), Vieweg,
  Braunschweig, 2006,  141--154.}\\ 

\vspace{-0.65cm}
\item{\emph{Fel'shtyn A., Indukaev F., Troitsky E.} {Twisted {B}urnside
theorem for two-step torsion-free nilpotent groups} // Preprint
MPIM2005-45, Max-Planck-Institut f{\"u}r Mathematik, 2005.}\\ 

\vspace{-0.65cm}
\item{\emph{Fel'shtyn A., Troitsky E., Vershik A.} {Twisted {B}urnside
theorem for type {II}${}_1$ groups: an example} // Math. Res. Lett. 2006.\textbf{13}, N 5--6 719--728.}\\

\vspace{-0.65cm}
\item{\emph{A.~A. Kirillov}, {Elements of the theory of representations},
  Springer-Verlag, Berlin Heidelberg New York, 1976.}\\ 

\vspace{-0.65cm}
\item{ \emph {A. O. Barut and R. Raczka}, Theory of Group Representations and Applications, Vol. 2 (PWN, 1977; Moscow,
Mir, 1980).}\\ 

\vspace{-0.65cm}
\item{\emph{A. O. Barut and R. Raczka}, Theory of Group Representations and Applications, Vol. 1 (PWN, 1977; Moscow,
Mir, 1980).}\\ 

\vspace{-0.65cm}
\item{\emph{Bernhard Amberg}. {Groups with maximum conditions} //  Pacific J.
Math.  1970.\textbf{32}, N.~1 9--19 (avaliable at
http://projecteuclid.org).}\\ 

\end{enumerate}
}

\end{document}